\documentclass[12pt]{article}

\usepackage{amssymb}
\usepackage{amsthm,amsfonts,amsmath,amssymb}

\newcommand{\eq}{\begin{equation}}
\newcommand{\en}{\end{equation}}

\newcommand{\prob}{\mathbb P}
\newcommand{\ex}{\mathbb E}

\newtheorem{theorem}{\large Theorem}
\newtheorem{proposition}[theorem] {\large Proposition}

\newtheorem{corollary}[theorem]{\large Corollary}

\begin{document}

\title{Recognising the Last Record of a Sequence}
\author {Alexander Gnedin\thanks{Insitute of Mathematics, Utrecht University, The Netherlands, gnedin@math.uu.nl }}
\maketitle

\begin{abstract}\noindent
We study the  best-choice problem for  processes which generalise the  
 process of records from   Poisson-paced 
i.i.d. observations. 
Under the assumption that  the
 observer knows  distribution of the process and the horizon,
we determine the optimal stopping policy and for a parametric family of 
problems also
derive an explicit formula for 
the maximum probability of recognising the last record.
\end{abstract}

\section{Introduction}
Maximising 
the probability of stopping at the extreme 
of a sequence 
of random marks
is the classical objective in sequential decision problems widely known as 
 the best-choice or `secretary' problems 
 \cite{BerGn, Obj, SamSurv}.
Problems of this kind can be formulated in terms of the embedded process of records, because the overall
extreme (e.g. minimum) is the {\it last record} observation.

\par In a basic version of the problem  introduced by Gilbert and Mosteller \cite[Section 3]{GM}
the  marks are sampled at  discrete times from the uniform distribution, and  the objective of the observer 
is to stop
at the  minimum among the first $n$ marks. 
The sequence of values of sequential minima, called  {\it lower records}, 
undergoes a {\it stick breaking}
process  $X_1, X_1X_2, X_1X_2X_3,\ldots$, where the $X_j$'s are independent copies of a prototypical random factor
$X$ whose distribution is uniform. 
Given the record values, the {\it durations} of records are independent,
and for $r$ a generic record value, the 
duration of a record with this value has geometric distribution with parameter $r$.
See \cite{Arnold, Goldie, Nevz} for these basic facts of the theory of random records.
The optimal policy in the 
stopping  problem of \cite{GM}
is rather complicated, as  it involves 
 a sequence of thresholds for which no closed-form expression is available,  and for the same reason
there is no explicit formula for the optimal probability.

\par According to another version of the problem, the marks are observed at epochs of a unit Poisson process,
and the goal  is to stop at the minimum mark before given horizon $T$
(see 
\cite{BC} and references therein). 
This problem allows much more explicit results: the optimal 
policy prescribes stopping at the first time the record process breaks through a hyperbolic 
boundary, 
and there is an explicit
formula for the optimal probability. 
The continuous time problem corresponds to the
 model sometimes called {\it Poisson-paced records\,}  \cite[Section 9]{Goldie}, the difference
with the discrete-time model is that
the duration of a record with value $r$ has exponential distribution 
with rate $r$.
For large $n$ and $T$ the discrete and continuous time  versions are close to
the same limiting form introduced in \cite{FI}, in particular the limiting
optimal probability is given by the  formula first obtained by Samuels \cite{SamFI} in discrete setting.
These and related results are reviewed in Section \ref{uniform}.

\par  In this paper we  consider the continuous-time problem of recognising the last record 
under a more general assumption that the occurences of records follow  a stick-breaking scheme,
with factor  $X$ having an arbitrary distribution on the unit interval. 
Models of this kind
appear in many contexts such as 
branching processes, search problems,  sequential packing problems  and random
partitions \cite{BarGn, Bertoin, Devroye, PTRF}.
Although we just
postulate the behaviour of records without any reference to some 
more rich observable process, the model in focus 
is related to one  concept of 
sequential extreme for sampling from certain  partially ordered spaces,
including spaces
${\mathbb R}^d$ with continuous product distributions. This connection is detailed in Section \ref{chain}.

\par  We will show that the optimal policy
is always of the same form as in the case of uniform $X$.
In one special case of  parametric family of beta distributions we 
express
the optimal probability in terms of the incomplete gamma function. 
In general, however,
it does not seem possible to write a  closed-form expression for the
stopping value. Still, we argue that under minor side conditions on the law of $X$, 
as $T\to\infty$, 
there exists a limiting value which may be interpreted as the 
optimal probability of recognising the last record in a stopping problem with infinitely many observations.
The famous best-choice probability benchmark $e^{-1}$ will show up as a sharp lower bound.

\section{The model}

We  shall model the occurences of records by means of
a nonincreasing right-continuous Markov process $R=(R_t,\,t\geq 0)$  with the following
type of behaviour: given the current state is $r>0$, the process jumps at rate $r$ to a new state $r X$, where
$X$ is a prototypical  random 
factor with a given distribution in the open interval $]0,1[\,$.
In the event $t$ is a jump instant of $R$ we say that a record occurs
at time  
$t$ and intepret $R_t$ as the  {\it weight} of the record. 
 The weights of  consequitive records decrease, while the sojourns
of $R$, which include the first record time and further durations of records,
 are stochastically increasing.
\par In more detail,
the weights of records undergo stick-breaking
$r_0X_1,r_0X_1X_2,\ldots$ , where  $X_j$'s
are independent replicas of  $X$ and $r_0=R_0$ is the initial state of $R$. 
The sequence of sojourns may be represented as
$$E_1/r_0, E_2/(r_0X_1),E_3/(r_0X_1X_2),\ldots$$
 where $E_j$'s are i.i.d. unit exponential variables, independent of the $X_j$'s.
Thus, conditionally  given the weights of records, the sojourns are independent 
exponential variables.

\par We are interested in the problem of maximising the probability of recognising 
the last record of $R$ before a given horizon $T$, by means of a nonanticipating policy (stopping time) 
adapted to the natural right-continuous filtration of $R$.
For  $\pi$ such a policy the efficiency is measured
by the probability that $\pi$ is a record time not exceeding $T$ and that 
no further record occurs before time $T$:
\eq\label{crit}
\prob(R_{\pi-}>R_\pi=R_T, \pi<T)=\ex\left[\exp\left\{-(T-\pi)R_\pi\right\}\,{\bf 1}(R_{\pi-}>R_\pi,\, \pi<T)\right],
\en
where the second expression involves the adapted probability 
of recognising the last record when the stopping occurs. 

\par In the terminology going back to  
Gilbert and Mosteller \cite{GM},
this stopping problem should be qualified as a problem with 
`full information', 
meaning that the observer  learns the weights of records and knows their distribution
exactly. 
Under `no-information' problem we understand the optimal stopping problem where
only policies based on record times are allowed.

\section{Chain records}\label{chain}

Sampling from arbitrary continuous distribution $F$ 
leads to the stick-breaking process for records with uniform $X$.
This is seen by defining the weight via 
$v\mapsto F(v)$ and by noting that this mapping
preserves the ranking and transforms a sample from $F$  into  a sequence of uniform variables.
In this section we discuss some extensions of this framework.

\par Sampling from certain discrete distributions also leads to stick-breaking process for records.
Define a distribution by allocating  the geometric masses $pq^{k-1}$  (where $p+q=1$ and $0<p<1$) 
at points of some  decreasing sequence $z_k$, $k=1,2,\ldots$.
Consider {\it strict} lower records in a sample from such distribution.
Define the weights 
by means of  the {\it left}-continuous
distribution function $v\mapsto F(v-)$.
If the first sample value is $z_k$, then the next observation is a record with probability $q^k$;
from this we see that the weights of records  
follow 
the stick-breaking scheme with factor
$$X=_d\sum_{k=1}^\infty pq^{k-1}\delta_{q^k}, $$
where $\delta_x$ is the Dirac mass at $x$ and $=_d$ denotes the equality in distribution.
\par Sampling from other distributions on reals is not consistent with the stick-breaking 
model for records. We will look now in higher dimensions.

\par Consider ${\mathbb R}^d$ endowed with some continuous product distribution $\mu$ and
the natural strict partial order $\prec$.
For a sample $V_1,V_2,\ldots$
 from $({\mathbb R}^d,\mu)$, 
we say that a {\it chain record} occurs at index $j$ if either $j=1$ or $j>1$ and 
$V_j$ is $\prec$-smaller than the last chain record in the sequence $V_1,\ldots,V_{j-1}$.
Define the {\it weight} of a chain record by means of  
the multivariate distribution function
$v\mapsto \mu\{u\in {\mathbb R}^d: u\prec v\}$.
The weights of chain records follow a stick-breaking process with the density  
$\prob(X\in{\rm d}x)/{\rm d}x=|\log x|^{d-1}/(d-1)!$
for the factor $X$.
Indeed,
the componentwise probability transform 
 establishes  isomorphism between the ordered probability
space $({\mathbb R}^d,\mu,\prec)$ and 
the unit cube $[0,1]^d$ with the Lebesgue measure, which implies that
the law of $X$ 
 is the same as the distribution 
of the product of $d$ independent uniform variables, whence the formula for the density.
\par Chain records in ${\mathbb R}^d$ were introduced in \cite{CR}.
Unlike other kinds of multidimensional records  surveyed in \cite{GR}, 
the chain records  cannot be regarded as  `generalised minima',
because 
permutations of $V_1,\ldots,V_{j-1}$ may destroy or create a chain record at index $j$.
The sequence of chain-record marks is a `greedy' decreasing chain in the partially ordered sequence 
of marks,  in the sense that element $V_j$ is joined to the chain each time when the monotonicity constraint
is not violated (as to be compared, e.g., with the longest chain  among the first $n$ marks).

\par The definition of chain record extends in an obvious way to sampling from an arbitrary Borel space ${\cal Z}$
endowed with a probability measure $\mu$ and a 
measurable strict partial order $\prec$. The weights are defined by 
means of the function $v\mapsto\mu(L_v)$, where $L_v=\{u\in {\cal Z}:u\prec v\}$ is the lower section of $\prec$
at $v\in {\cal Z}$.
Call the space $({\cal Z},\mu,\prec)$ {\it lower-homogeneous} if (i) $\mu(L_v)>0$
 for $\mu$-almost all points $v\in {\cal Z}$, and (ii) 
the lower section $L_v$ with conditional measure $\mu(\cdot)/\mu(L_v)$
is isomorphic, as a partially ordered probability space, to the whole space
$({\cal Z},\mu,\prec)$. Since all $L_v$'s are in this sense the same,
the weights of chain records in a sample from a lower-homogeneous space undergo a stick-breaking
with the factor $X=_d\mu\{u\in {\cal Z}:u\prec V\}$ where $V$ has distribution $\mu$.

\par It is easily seen that  $[0,1]^d$ with uniform distribution 
is a lower-homogeneous space,
hence  this is  true also for  ${\mathbb R}^d$ with continuous product distribution.
Another example is 
the {\it interval space} which has intervals $]a,b[\subset [0,1]$ as elements,
the partial ordering $\prec$ defined by inclusion, and a measure  $\mu({\rm d}a{\rm d}b)=\alpha(\alpha-1)(b-a)^\alpha{\rm d}a\,{\rm d}b$ (with
parameter
$\alpha>1$); in this case
$\prob(X\in {\rm d}x)/{\rm d}x=(\alpha-1)(x^{-1/\alpha}-1)$.
Although both examples are instances of Bollob{\'a}s-Brightwell {\it box-spaces} \cite{BB}
(which have all
intervals
$\{u: v\prec u\prec w\}$ for $v\prec w$
 isomorphic to the whole space  and not only $L_v$'s),
there are many other lower-homogeneous spaces that are not box-spaces.
By the transitivity of partial order, 
the distribution  of   
$X$ appearing in this way
must satisfy the inequality $\prob(X\leq x)\geq x\,,~x\in[0,1]$.

\section{Stopping the embedded Markov chain}
\label{embed}
 A fundamental property of the process $R$ is {\it self-similarity}:
for each $r>0$,
the law of $(R_t,~t\geq 0)$ with $R_0=r$ is identical to the law of $(r R_{rt},~t\geq 0)$ given $R_0=1$.
This implies that, when the law of $X$ is fixed, the `size' of the problem is determined by a single
parameter $r_0T$.

\par Self-similarity is  a clue to derive the optimal policy.
If stopping has not occur before and including time $t$ and if the 
current state is $R_t=r$,
then the conditional optimal stopping problem is equivalent to the unconditional  problem with
initial state $1$ and horizon $(T-t)r$.
 This motivates associating with $R$ (with fixed parameters $r_0,T$ and the law for $X$)
another decreasing Markov process $B=(B_t,\,t\geq 0)$,
$$B_t=(T-t)_+R_t\,,~~~t\geq 0,$$
with the initial state $B_0=r_0T$ and the absorbing terminal state $0$. 
Obviously, it is sufficient to consider
the policies adapted to $B$, with understanding
that  the last record before $T$ corresponds to
the last jump of $B$ before absorption at $0$.
The sequence of locations visited by $B$ at  the record times  is a discrete-time homogeneous Markov chain which
follows the transition scheme
$s\mapsto (s-E)_+X$, $s\geq 0$, where $E$ is a rate-$1$ exponential variable independent of $X$.

\par 
In terms of the embedded chain the optimal policy is determined in a standard way, 
by comparing two kinds of risk.
If the current state of $B$ is $s$, the probability that no further records occur
is $p_0(s)=e^{-s}$. On the other hand, the probability that exactly one record will occur is
$$p_1(s)=\int_0^s e^{-t} \,\ex\,[p_0((s-t)X)]\,{\rm d}t=e^{-s}\, \ex\left[{e^{s(1-X)}-1\over 1-X}\right].$$
Inspecting two extremes $s=0$ and $\infty$ and exploiting  monotonicity,  we see that the equation
\eq\label{thr}
\ex\left[{e^{s(1-X)}-1\over 1-X}\right]=1
\en
has a unique positive solution $s_*$. Because 
$$p_0(s)<p_1(s)\Longleftrightarrow s>s_*\,,$$
and because $B$ has decreasing paths (until getting absorbed)
we are in the familiar  {\it monotone case} of optimal stopping, hence 
the optimal policy stops at the first jump of  $B$ within the region $[0,s_*]$.
Translating this back in terms of $R$ we see that it is optimal to stop at the
first record time when the condition
$(T-t)R_t\leq s_*$
is satisfied. In particular, if $T r_0\leq s_*$ it is optimal to stop at the very first record.

\par More generally, for $s>0$ we denote $\pi_s$ the policy  which prescribes stopping at the first record time
when $(T-t)R_t\leq s$ holds.
Summarising  the above discussion  we conclude:

\begin{proposition}
The policy $\pi_{s_*}$ with $s_*$ satisfying {\rm (\ref{thr})} 
is optimal.
\end{proposition}

\vskip0.5cm

Assuming $R_0=1$,
let $v(T,s)$ be
the value of the policy $\pi_s$, i.e.
 the probability that exactly one record before $T$ satisfies $(T-t)r\leq s$.
(By self-similarity  the case of arbitrary $R_0=r_0$ can be reduced to that.)
Obviously, 
\eq\label{init}
v(T,s)=p_1(T)\,,~~~{\rm  for~~} T<s\,.
\en
The first-record decomposition  readily yields
an integral equation
$$
v(T,s)=\int_0^{T}\,\ex\left[v(X(T-t),s) \,
{\bf 1}((T-t)X>s)+e^{-(T-t)X}\,{\bf 1}((T-t)X\leq s)
\right] e^{-t}{\rm d}t\,,
$$
which for $s=s_*$ is the familiar dynamic programming equation for the optimal value.
In the differential form this becomes
\eq\label{diff}
\partial_T \,v(T,s)=-v(T,s)+
\ex\left[v(TX,s)\,{\bf 1}(TX>s)\right]+
\ex\left[e^{-TX}\,{\bf 1}(TX\leq s)\right].
\en

\par The equation (\ref{diff}) is of delayed type, which only in exceptional cases admits a closed-form
solution. For instance, when the distribution of $X$ is $\delta_x$, the solution is a piecewise-analytical
function which  should be computed recursively in the intervals $T\in [s/x^{k-1},s/x^k]$ for $k=1,2,\ldots$,
starting from $[0,s/x]$ where $v(T,s)=p_1(s)$ holds.
\par The collection of sites which $B$ visits at record times is not a Poisson
process, since otherwise $v(T,s)$ were constant in $T$ for $T>s$.
It is therefore surprising 
that the maximum of $v(T,s)$  in $s$ is attained at the same point $s_*$,
for all $T>s_*$.

\section{The lower bound}

Suppose for a while that the law of $X$ is $\delta_1$. In this case (\ref{diff}) is easily solved as
$v(T,s)=(T\wedge s)\,e^{-(T\wedge s)}$. Thus $s_*=1$ and for $T\geq 1$ the optimal probability is 
$v(T,s_*)=e^{-1}$, which  also coincides with the maximum of  $p_1(s)=se^{-s}$.
To bring this conclusion into the familiar `no-information' framework
 note that $R_t\equiv 1$, 
hence there is no updating of record weights. For the same reason, the record times
are the epochs of a unit Poisson  process, hence the 
 stopping 
problem amounts  to recognising the last Poisson epoch on $[0,T]$, which is 
 the `no-information'  problem
for Poisson process   due to Browne \cite{Browne}.
A characteristic feature of this case is that $v(T,s)$
is constant in $T$ for $T>s$ (see the last remark in Section \ref{embed}).

\par We show next that the familiar
 benchmark $e^{-1}=0.367\ldots$  yields a universal lower bound 
in our model.

\begin{proposition} For 
every distribution of $X$
the optimal probability satisfies
$v(T,s_*)>e^{-s_*}$ for $T>s^*$.
Above that, $s_*<1$ hence 
$$v(T,s_*)>e^{-1}\,~~{\rm for ~\,}T>1,$$
and this bound is sharp.
\end{proposition}
\begin{proof} Suppose $r_0T>s_*$. The process $B$ can enter $[s_*,0]$ by either continuously 
drifting down or jumping down through $s_*$.
In the first case the conditional probability of sucess with $\pi_*$ is $p_1(s_*)=p_0(s_*)=e^{-s_*}$.
 In the second case this probability is
$\ex[e^{-S}]>e^{-s_*}$, with some random $S<s_*$, because 
$\pi_{s_*}$ will stop. 
The estimate readily follows.

\par Applying the inequality $e^{1-x}>1+(1-x)$ for $0<x<1$, we see that
the left-hand side of (\ref{thr}) is larger than $1$ for $s=1$,
therefore the root satisfies $s_*<1$. The bound $e^{-1}$ is approached by letting  
the law of $X$ to approach $\delta_1$.
\end{proof}

\noindent
 The same argument yields a more general inequality $v(T,s)>\min(p_0(s),p_1(s))$ for $T>s$, where the right side assumes the 
largest value at $s=s_*$.

\par For $X$ uniform $s_*=0.804\ldots$ and the lower bound is $e^{-s_*}=0.447\ldots$,
while for $X$ with density $|\log x|$ these are $0.743\ldots$ and $0.475\ldots$.

\section{Entrance from the infinity}
\label{entr}

For asymptotic considerations we shall vary the initial state and denote $\prob_{r}$ the law of $R$
with $R_0=r$.
Assume that $\ex\,|\log X|<\infty$  and that the distribution of $X$ is not supported by a geometric progression
(note that these are precisely the conditions for applicability of the renewal theorem \cite{Feller} to $-\log X$).
Let  $$f(\lambda)=\ex[X^\lambda]$$
be the Mellin transform of $X$. Clearly, $-f'(0)=\ex\,|\log X|$. 
Adapting \cite[Theorem 1]{Entr} we have:

\begin{proposition}\label{entrance} Under the above assumptions, as $r\to \infty$, the law $\prob_r$  
has a weak limit $\prob_\infty$ characterised by $R_t=_d Y/t\,,\,t>0,$ where $Y$ 
is a random variable uniquely determined by its 
 moments 
\eq\label{mom}
\ex \,[Y^k]={1\over -f'(0)}\prod_{j=1}^{k-1} {j\over 1-f(j)}\,\,,~~k=1,2,\ldots
\en
\end{proposition}

\begin{corollary}
Under these circumstances
there exists a limit $v(\infty,s_*)=\lim_{T\to\infty} v(T,s_*)$ which is the maximum probability
of recognising the last record for the process $(R_t,\,t\leq 1)$ under  $\prob_\infty$.
\end{corollary}
\begin{proof} This follows from the form of the optimal policy and the fact that 
the point process of sites visited by $B$ at record times has a weak limit as $B_0\to\infty$.
\end{proof}
\noindent
The law of $Y$ determined by (\ref{mom}) can be considered as a kind of extreme-value distribution.
For instance, $Y$ is exponential for $X$ uniform, while $Y$ is distributed like the product of independent uniform
and exponential variables for $X$ with density $|\log x|$. 

\par Denoting $\tau_1,\tau_2,w_1,w_2$ the times and weights of the last record and the record before the last,
the performance of $\pi_s$ in the infinite problem can be written as
\eq\label{vinfty}
v(\infty,s)=\prob_\infty((1-\tau_1)\rho_1<s<(1-\tau_2)\rho_2)
=\prob_\infty((1-\tau_1)\rho_1<s)-\prob_\infty((1-\tau_2)\rho_2<s).
\en
In principle, the moments (\ref{mom}) determine the distribution of these  variables,
for instance
$$\prob_\infty(\tau_1<t)=\ex\left[e^{Y(1-1/t)}\right],$$
but it seems impossible to use this for writing $v(\infty,s)$ in some explicit form.

\section{The beta case}

We proceed with more concrete computations under the assumption that the distribution of $X$ is beta$(\theta,1)$,
with the density
$$\prob(X\in {\rm d}x)/{\rm d}x=\theta x^{\theta-1}\,,~~~x\in [0,1],$$
where  $\theta$ is a positive parameter. 
The instance $\theta=1$ corresponds to the uniform distribution.
This class of stick-breaking processes has a feature that
 under $\prob_\infty$ both the range of $R$ and the point process of record times
are Poisson point processes with intensity measure $\theta\,{\rm d}z/z$, $z>0$.
The law of $R_1$  under $\prob_\infty$ is a gamma distribution.
\par The integral 
$$p_1(s)=\int_0^1 {e^{-sx}-e^{-s}\over 1-x}\,\theta x^{\theta-1}\,{\rm d}x\,$$
does not simplify, hence it should be included in the final formula for the optimal probability as it is.

\subsection{Computing the value}

For $T>s$  a substitution translates (\ref{diff}) into
\begin{eqnarray*}
T^\theta\,\partial_T\, v(T,s)=-T^\theta v(T,s)+\int_s^T v(t,s)\theta t^{\theta-1}{\rm d}t
+\int_0^s e^{-t}\theta t^{\theta-1}{\rm d}t\,.
\end{eqnarray*}
Differentiating in $T$ and simplifying we are lead to 
\eq\label{g}
Tg''+(T+\theta)g'=0
\en
for $g(T)=v(T,s)$. Solving this and taking into account the boundary condition at $T=s$
yields
\eq\label{val}
v(T,s)=\Gamma(-\theta+1, s,T)e^s s^\theta p_1'(s)+p_1(s)\,,~~{\rm for~~} T>s\,,
\en
where
$$\Gamma(a,b,c)=\int_b^c e^{-t}t^{a-1}{\rm d}t$$
denotes the incomplete gamma function,
and
$$p_1'(s)=-p_1(s)+ {\theta\over s^\theta}\,\Gamma(\theta,0,s).$$

\par For the optimal $s_*$ using $p_1(s_*)=e^{-s_*}$ we obtain  from 
(\ref{val}) 
\eq\label{value}
v(T,s_*)=\Gamma(-\theta+1,s_*,T)[-s_*^\theta+e^{s_*}\theta\Gamma(\theta,0,s_*)]+e^{-s_*}\,,
\en
which is the optimal probability
of stopping at the last record. 
The formula is valid for $T\geq s^*$.
The optimal probability $v(\infty,s_*)$ in the limit problem
is just obtained 
taking $T=\infty$  in the integral in (\ref{value}), which reads as
a generalised exponential integral function
$$\Gamma(-\theta+1,s_*,\infty)=\int_{s_*}^\infty {e^{-t}\over t^\theta}\,{\rm d}t\,.$$

\par The following table shows some numerical values of this probability computed with a help of 
{\tt Mathematica}.

\begin{eqnarray*}
\begin{array}{clllllll}
\theta &0.1&0.25&0.5&1&2&5&20\\
s_* & 0.709& 0.731& 0.760& 0.804& 0.857& 0.922& 0.976\\
v(\infty,s_*)& 0.913& 0.814&0.703 & 0.580& 0.481& 0.410& 0.377 
\end{array}
\end{eqnarray*}
The  data suggest   to examine the extreme values of the parameter $\theta$.

 \par As $\theta\to\infty$ the beta distribution approaches $\delta_1$, hence $s_*\to 1$ and
$v(\infty,s_*)\to e^{-1}$. Thus the beta family may be interpreted as a bridge between the
`full-information' problem ($\theta=1$) and the `no-information' problem ($\theta=\infty$).

\par Note that, for arbitrary $\theta>0$, in consequence of the Poisson character of the record times 
under $\prob_\infty$,
the time-threshold policy $\pi=\min\{t>T/e,\,R_{t}>R_{t-}\}$  yields the limit probability
of success equal $e^{-1}$ for $T\to\infty$.

\par As $\theta\to 0$ the beta distribution approaches $\delta_0$.
In this regime the optimal $s_*$ approaches $\log 2$.
Selecting $T_0$ sufficiently large to secure occurence of at least one record with probability at least $1-\epsilon$,
and then sending $\theta$ to $0$, we will have $v(T_0,s_*)\geq 1-\epsilon$,
because with high probability exactly one record  occurs before horizon $T$. 
For $T>T_0$ (\ref{value}) implies $v(T,s_*)>v(T_0,s_*)$, therefore
the trivial upper bound $v(\infty,s_*)<1$ is sharp as the law of $X$ varies.

\subsection{A smooth fit}

With the explicit formula (\ref{val}) in hand we can alternatively 
characterise $s_*$ as the maximiser of  $v(T,s)$ in $s$. Equating
$\partial_s\,v(T,s)$ to $0$ 
we see then that $s_*$ is a root of the equation
\eq\label{r1}
s p_1''(s)+(s+\theta)p_1'(s)=0\,,
\en
which is, in fact, equivalent to
$p_0(s)=p_1(s)$ due to the  identity
$$sp_1''(s)+(s+\theta)p'(s)=\theta(p_0(s)-p_1(s)).$$

\par Comparing (\ref{r1}) with (\ref{g}) shows that two branches of $v(T,s_*)$,
for $T\leq s_*$ and $T\geq s_*$, match at $T=s_*$ together with  {\it two} derivatives.
This degree of smoothness is characteristic for $s=s_*$, as is also seen by the following argument. 
Write the probability of success with policy $\pi_s$ as
$$v(T,s)=\int_0^s p_1'(t){\rm d}t+\int_s^T \partial_T \,v(t,s){\rm d}t\,,$$
and note that
the optimisation of $s$ amounts to finding a `switch' which maximises the sum of integrals.
Inspecting  
the monotonicity properties of the integrands 
$$p_1'(t)={e^{-t}\over t^\theta}(p_1'(t)t^\theta e^{t})~~~{\rm and~~}~\partial_T\,
v(t,s)={e^{-t}\over t^\theta}(p_1'(s)s^\theta e^{s})$$
shows that the maximum is achieved if they are
tangential at the switching location, 
which is precisely the condition (\ref{r1}).
Thus $s_*$ 
is indeed the only value of $s$ such that
$\partial_T\,v(T,s)$ has no break at $T=s$.

\subsection{The uniform case}
\label{uniform}

For completeness we bring together known formulas for  the case of 
uniform factor $X$.
\par The solution to 
$$\int_0^s {e^t-1\over t}{\rm d}t=1$$
has the approximate value $s_* =0.804\ldots$
The limit probability
$$v(\infty,s_*)=(e^{s_*}-s_*-1)\int_{s_*}^\infty {e^{-t}\over t}\,{\rm d}t+e^{-s_*}=0.580\ldots\,,$$
was obtained first numerically in \cite{GM} by interpolation from discrete-time problems, 
derived
in \cite{SamFI} from (\ref{vinfty}),
and shown \cite{BerGn} by some series computations with the Poisson process.
The  analogous formula for $v(s,T)$ with  finite $T$ appeared in \cite{BC}.

\par The process of records under $\prob_\infty$  corresponds to the set of $\prec$-minimal atoms of a unit-rate
 Poisson point process on ${\mathbb R}_+^2$ 
(recall that an atom is $\prec$-minimal if there are no other Poisson atoms south-west of it).
The unique  properties of the planar Poisson process allow more delicate computations.
Under $\prob_\infty$ 
the density of $\pi_s$ is \cite{BC}
\begin{eqnarray*}
\prob(\pi_s\in {\rm d}t)/{\rm d}t=
{t-1\over t}(e^{-ts}-e^{-ts/(1-t)})+s\Gamma\left(0,st,{st\over 1-t}\right)+1-e^{-st}\,,~~t\in [0,1],
\end{eqnarray*}
This integrates to some number less than $1$  because  with positive 
probability $\pi_s$ does not stop at all (this probability is aproximately $0.1995\ldots$ for the optimal policy $\pi_{s_*}$).
The optimal probability can be also represented as the integral
$$v(\infty,s_*)=\int_0^1 w(t){\rm d}t\,,$$
with $w(t)$  the
 {\it winning rate}, equal to the chance that $\pi_{s_*}$ stops correctly
in time ${\rm d}t$.
The graph of $w$ was sketched 
long ago  \cite[Figure 3]{GM}, and
the following explicit formula for the winning rate is a recent result \cite{GnMir}:
\begin{eqnarray*}
w(t)=-e^{-s_*}+{e^{-s_* t}-e^{-s_*t/(1-t)}\over t}+{e^{-s_*t}-te^{-s_*}\over 1-t}
+\\
{s_*\over 1-t}\left[\Gamma\left(0,s_*,{s_*\over 1-t}\right)-\Gamma\left(0,s_*t,{s_* t\over 1-t}\right)\right],
\end{eqnarray*}
(the boundary values $w(0)=1-e^{-s_*}, ~w(1)=e^{-s_*}$ were indicated in  \cite{GM}).

\section{Concluding remarks}

A discrete-time version of the problem with fixed horizon $n$
is associated with a process analogous to $R$ but with geometric durations of records.
The optimal policy is known only for uniform $X$.
 Moreover, it is not clear if the monotone 
case of optimal stopping applies for the general  distribution of $X$.
Using techniques from
\cite{KR}
one can show that
under the assumptions of Section \ref{entr}  the discrete-time problem
can be aproximated, for $n$ large, 
by the limiting problem with continuous time,
hence the policy
`stop at index $j$ if a record occurs with weight
$r$ satisfying $r(n-j)\leq s_*$' is
asymptotically optimal.


\par It would be also interesting to evaluate suboptimal policies like `stop at the first record
with weight below given $w$' or `stop at the first record that occurs after a given time $t_0$'.
This is not so easy in general since such policies are not adapted to  $B$.

\end{document}